\documentclass[12pt]{amsart}

\usepackage{latexsym}
\usepackage{amssymb}

\newcommand{\excise}[1]{}

\newtheorem{thm}{Theorem}[section]
\newtheorem{lemma}[thm]{Lemma}

\newtheorem{cor}[thm]{Corollary}
\newtheorem{prop}[thm]{Proposition}

\newtheorem{question}[thm]{Question}
\newtheorem{problem}[thm]{Problem}
\newtheorem{por}[thm]{Porism}

\newtheorem{Example}[thm]{Example}
\newtheorem{Remark}[thm]{Remark}
\newtheorem{Alg}[thm]{Algorithm}
\newtheorem{Defn}[thm]{Definition}

\newenvironment{example}{\begin{Example}\rm}
                {\mbox{}~\hfill$\square$\end{Example}}
\newenvironment{remark}{\begin{Remark}\rm}
                {\mbox{}~\hfill$\square$\end{Remark}}

\newenvironment{defn}{\begin{Defn}\rm}
                {\mbox{}~\hfill$\square$\end{Defn}}

\newenvironment{proofof}[1]{\begin{trivlist}\item {\bf
        Proof of {#1}.\,}}{\mbox{}\hfill$\square$\end{trivlist}}

	{
	 \setcounter{separated}{\value{equation}}
	 \setcounter{equation}{0}
	 \begin{equation}
	}
	{\end{equation}%
	 \setcounter{equation}{\value{separated}}%
	}

\newenvironment{numbered}%
        {\begin{list}
                {\noindent\makebox[0mm][r]{\arabic{enumi}.}}
                {\leftmargin=5.5ex \usecounter{enumi}}
        }
        {\end{list}}

\newenvironment{romanlist}%
        {\begin{list}
                {\noindent\makebox[0mm][r]{(\roman{enumi})}}
                {\leftmargin=5.5ex \usecounter{enumi}}
        }
        {\end{list}}

\newcounter{separated}

\def\bem#1{\textbf{#1}}

\def\<{\langle}
\def\>{\rangle}
\def\0{{\mathbf 0}}
\def\1{{\mathbf 1}}

\def\cD{{\mathcal D}}

\def\FF{{\mathcal F}}
\def\GG{{\mathcal G}}
\def\HH{{\widetilde H}{}}

\def\KK{{\mathcal K}}

\def\NN{{\mathbb N}}

\def\SS{{\mathfrak S}}

\def\ZZ{{\mathbb Z}}

\def\bb{{\mathbf b}}

\def\kk{{\mathbf k}}

\def\xx{{\mathbf x}}
\def\yy{{\mathbf y}}
\def\zz{{\mathbf z}}

\def\th{{\rm th}}
\def\abs{{\rm abs}}

\def\del{{\rm del}}

\def\link{{\rm link}}

\def\length{{\rm length}}
\def\repnum{{\rm repnum}}

\def\K{$K$}

\def\from{\leftarrow}

\def\minus{\smallsetminus}

\def\goesto{\rightsquigarrow}
\def\implies{\Rightarrow}
\def\nothing{\varnothing}

\def\ol#1{{\overline {#1}}}

\def\dem#1{{\ol \partial_{#1}}}

\def\twoline#1#2{\aoverb{\scriptstyle {#1}}{\scriptstyle {#2}}}

\newcommand{\aoverb}[2]{{\genfrac{}{}{0pt}{1}{#1}{#2}}}

\font\co=circle10

\def\jr{\smash{\raise2pt\hbox{\co \rlap{\rlap{\char'005} \char'007}}
               \raise6pt\hbox{\rlap{\vrule height6.5pt}}
               \raise2pt\hbox{\rlap{\hskip4pt \vrule height0.4pt depth0pt
                width7.7pt}}}}
\def\je{\smash{\raise2pt\hbox{\co \rlap{\rlap{\char'005}
                \phantom{\char'007}}}\raise6pt\hbox{\rlap{\vrule height6pt}}}}
\def\+{\smash{\lower2pt\hbox{\rlap{\vrule height14pt}}
                \raise2pt\hbox{\rlap{\hskip-3pt \vrule height.4pt depth0pt
                width14.7pt}}}}

\def\textcross{\ \smash{\lower4pt\hbox{\rlap{\hskip4.15pt\vrule height14pt}}
                \raise2.8pt\hbox{\rlap{\hskip-3pt \vrule height.4pt depth0pt
                width14.7pt}}}\hskip12.7pt}
\def\textelbow{\ \hskip.1pt\smash{\raise2.8pt%
                \hbox{\co \hskip 4.15pt\rlap{\rlap{\char'005} \char'007}
                \lower6.8pt\rlap{\vrule height3.5pt}
                \raise3.6pt\rlap{\vrule height3.5pt}}
                \raise2.8pt\hbox{%
                  \rlap{\hskip-7.15pt \vrule height.4pt depth0pt width3.5pt}%
                  \rlap{\hskip4.05pt \vrule height.4pt depth0pt width3.5pt}}}
                \hskip8.7pt}

\begin{document}

\title{Subword complexes in Coxeter groups}
\author{Allen Knutson}
\thanks{AK was partly supported by the Clay Mathematics Institute, Sloan
        Foundation, and NSF}
\address{Mathematics Department\\ UC Berkeley\\ Berkeley, California}
\email{allenk@math.berkeley.edu}
\author{Ezra Miller}
\thanks{EM was supported by the NSF}
\address{Mathematical Sciences Research Institute\\ Berkeley, California}
\email{emiller@msri.org}
\date{11 March 2003}

\begin{abstract}
\noindent
Let $(\Pi,\Sigma)$ be a Coxeter system.  An ordered list of elements in
$\Sigma$ and an element in $\Pi$ determine a {\em subword complex}, as
introduced in \cite{grobGeom}.  Subword complexes are demonstrated here
to be homeomorphic to balls or spheres, and their Hilbert series are
shown to reflect combinatorial properties of reduced expressions in
Coxeter groups.  Two formulae for double Grothendieck polynomials, one
of which appeared in~\cite{FKgrothYangBax}, are recovered in the context
of simplicial topology for subword complexes.  Some open questions
related to subword complexes are presented.
\end{abstract}

\maketitle

{}


\section{Introduction}

\label{sec:intro}

We introduced subword complexes
in \cite{grobGeom}
to elucidate the combinatorics of determinantal ideals and Schubert
polynomials.  In retrospect, however, 
they raise basic questions about the nature of reduced expressions in
arbitrary Coxeter groups.  For instance, given a fixed word~--- that is,
an ordered list $Q$ of simple reflections in a Coxeter group~$\Pi$~---
what can be said about the set of all of its reduced subwords?
In particular, given a fixed element $\pi \in \Pi$, what structure
belongs to the set of subwords of~$Q$ that are reduced expressions
for~$\pi$?

The exchange axiom in~$\Pi$ answers this last question when $Q$ is a list
of $1+\length(\pi)$ simple reflections (see Lemma~\ref{lemma:T}).
The general answer, when $Q$ and~$\pi$ are arbitrary, lies in properties
of the subword complex $\Delta(Q,\pi)$, whose facets correspond (by
definition) to reduced subwords of~$Q$ having product~$\pi$.

Both the topology and combinatorics of $\Delta(Q,\pi)$ are governed by
the exchange axiom in a strong sense.
Our first main result, Theorem~\ref{thm:ball}, says that subword
complexes are homeomorphic to balls or spheres.  The proof uses the fact
that subword complexes are shellable, which was demonstrated in
\cite{grobGeom} by exhibiting an explicit vertex decomposition.  The
lurking exchange axiom surfaces here as the transition between adjacent
facets across the codimension~$1$ face joining them.

Given their
topological simplicity, the invariants of subword complexes necessarily
derive from more refined combinatorial data, namely the links of all
faces.  Therefore we focus on homological aspects of Stanley--Reisner
theory in
Section~\ref{sec:hilb}, where we calculate the Hilbert series of face
rings of subword complexes in Theorems~\ref{thm:hilb}
and~\ref{thm:hilb'}.  The exchange axiom here gives rise to the
criterion for a face to lie in the boundary of a subword complex.

The structure theorem for subword complexes, Theorem~\ref{thm:ball}, is
reminiscent of fundamental
results for Bruhat as well as weak orders, and also for finite
distributive lattices.
The topology of subword complexes looks similar to that of order
complexes of intervals in the Bruhat order studied by Bj\"orner and
Wachs \cite{BWbruhCoxetShel}, even though there seems to be little
direct connection to subword complexes.  Indeed, the simplicial
complexes in \cite{BWbruhCoxetShel} are by definition independent of the
reduced expressions involved, although the lexicographic shellings there
depend on such choices.  In contrast, Bj\"orner \cite[Section~3]{Bjo84}
concerns himself with intervals in weak orders, where the reduced
expressions involved form the substance of the simplicial complexes, as
they do for subword complexes.  Bj\"orner proves that intervals in the
weak order are homotopy-equivalent to balls or spheres.  Our results are
geometrically somewhat stronger, in that we prove not just homotopy
equivalence, but homeomorphism.

The comparison between subword complexes and 
order complexes of intervals in the weak order occurs most clearly when
the word $Q$ contains as a subword every reduced expression for~$\pi$.
In essence, the reduced expressions for~$\pi$ must be repeated often
enough, and in enough locations inside~$Q$, to make $\Delta(Q,\pi)$
homeomorphic to a manifold, whereas the set of reduced expressions
for~$\pi$~--- without repeats~--- only achieves homotopy equivalence.
Some open questions in Section~\ref{sec:open} are relevant here.

The plan of the paper is as follows.  We review in
Section~\ref{sec:subword} the shelling construction and its proof from
\cite{grobGeom}, along with related definitions.  Next we prove the
structure theorem in Section~\ref{sec:balls}.  Section~\ref{sec:hilb}
contains the Hilbert series calculation, which requires a review of
Hochster's formula and the Alexander inversion formula \cite{Mil3}, the
latter expressing the simple relation between the Hilbert series of a
squarefree monomial ideal and that of its Alexander dual.  In
Section~\ref{sec:groth}, we apply our Hilbert series formula in the
context of
symmetric groups from \cite{grobGeom} to deduce two formulae for
Grothendieck polynomials, one of which is due originally to Fomin and
Kirillov~\cite{FKgrothYangBax}.
Finally, we present some open problems in Section~\ref{sec:open}.

We note that the particular subword complexes in Example~\ref{ex:square}
were in fact the special cases that originally led us to define subword
complexes in general.  These special cases appear as the initial schemes
of certain types of determinantal varieties (`matrix Schubert
varieties').  The shellability for subword complexes proved in
\cite{grobGeom} and reviewed in Section~\ref{sec:subword} allowed us to
give in \cite{grobGeom} an independent proof of Cohen--Macaulayness for
matrix Schubert varieties \cite{FulDegLoc}, and therefore also for
ordinary Schubert varieties in the flag manifold \cite{ramanathanCM}.

\subsection*{Acknowledgements}
The authors are grateful to Sara Billey, Francesco Brenti, and Richard
Stanley for helpful discussions.  Cristian Lenart kindly suggested we
state the alternate formula for Grothendieck polynomials in
Corollary~\ref{cor:links'}, which prompted us to include also
Theorem~\ref{thm:hilb'}.

\section{Subword complexes}

\label{sec:subword}

We deal with an arbitrary Coxeter system $(\Pi,\Sigma)$ consisting of a
Coxeter group $\Pi$ and a set $\Sigma$ of simple reflections minimally
generating~$\Pi$.  See \cite{HumphCoxGrps} for background.
In Section~\ref{sec:groth}, we shall be particularly interested in an
application where where $\Pi = S_n$ is the symmetric group, and $\Sigma$
consists of the adjacent transpositions $s_1,\ldots,s_{n-1}$, where $s_i$
switches $i$ and $i+1$.
Here is our main definition, copied from
\cite[Definition~1.8.1]{grobGeom}.

\begin{defn} \label{defn:subword}
A \bem{word} of size~$m$ is an ordered sequence $Q = (\sigma_1, \ldots,
\sigma_m)$ of elements of~$\Sigma$.  An ordered subsequence $P$ of~$Q$
is called a \bem{subword} of~$Q$.
\begin{numbered}
\item
$P$ \bem{represents} $\pi \in \Pi$ if the ordered product of the simple
reflections in $P$ is a reduced decomposition for $\pi$.
\item
$P$ \bem{contains} $\pi \in \Pi$ if some subsequence of $P$ represents
$\pi$.
\end{numbered}
The \bem{subword complex} $\Delta(Q,\pi)$ is the set of subwords
$Q\minus P$ whose complements $P$ contain $\pi$.
\end{defn}
In other words, if~$Q\minus D$~is a facet of the subword
complex~$\Delta(Q,\pi)$, then the reflections in~$D$ give a reduced
expression for~$\pi$.  Note that subwords of~$Q$ come with their
embeddings into~$Q$, so two subwords $P$ and~$P'$ involving reflections
at different positions in~$Q$ are unequal, even if the sequences of
reflections in $P$ and~$P'$ are~equal.

Often we write~$Q$ as a string without parentheses or commas, and abuse
notation by saying that $Q$ is a word in~$\Pi$, without explicit
reference to~$\Sigma$.
Note that $Q$ need not itself be a reduced expression.  The following
lemma is immediate from the definitions and the fact that all reduced
expressions for $\pi \in \Pi$ have the same length.

\begin{lemma} \label{lemma:represents}
$\Delta(Q,\pi)$ is a pure simplicial complex whose facets are the
subwords $Q \minus P$ such that $P \subseteq Q$ represents
$\pi$.\hfill$\square$
\end{lemma}

\begin{example} \label{ex:pentagon}
Consider the subword complex $\Delta = \Delta(s_3s_2s_3s_2s_3,1432)$ in
$\Pi = S_4$.  This $\pi = 1432$ has two reduced expressions, namely
$s_3s_2s_3$ and $s_2s_3s_2$.  Labeling the vertices of a pentagon with
the reflections in~$Q = s_3s_2s_3s_2s_3$ (in cyclic order), we find that
the facets of $\Delta$ are the pairs of adjacent vertices.  Therefore
$\Delta$ is the pentagon.%
\end{example}

\begin{defn} \label{defn:vertex}
Let $\Delta$ be a simplicial complex and $F \in \Delta$ a face.
\begin{numbered}
\item
The \bem{deletion} of $F$ from $\Delta$ is $\del(F,\Delta) = \{G \in
\Delta \mid G \cap F = \nothing\}$.
\item
The \bem{link} of $F$ in $\Delta$ is $\link(F,\Delta) = \{G \in \Delta
\mid G \cap F = \nothing$ and $G \cup F \in \Delta\}$.
\end{numbered}
$\Delta$ is \bem{vertex-decomposable} if $\Delta$ is pure and either
(1)~$\Delta = \{\nothing\}$, or (2)~for some vertex $v \in \Delta$, both
$\del(v,\Delta)$ and $\link(v,\Delta)$ are vertex-decomposable.  A
\bem{shelling} of $\Delta$ is an ordered list $F_1, F_2, \ldots, F_t$ of
its facets such that $\bigcup_{j < i} \hat F_j \cap \hat F_i$ is a
subcomplex generated by codimension~1 faces of~$F_i$ for each $i \leq t$,
where~$\hat F$ denotes the set of faces of~$F$.  We say that $\Delta$ is
\bem{shellable} if it is pure and has a shelling.
\end{defn}

Vertex-decomposability can be seen as a sort of universal property.
Indeed, suppose that $\FF$ is a family of pure simplicial complexes in
which every nonempty complex $\Delta \in \FF$ has a vertex whose link
and deletion both lie in~$\FF$.
Then $\FF$ consists of vertex-decomposable complexes.  The set of
vertex-decomposable complexes is the largest (hence universal) such
family.

In the above definition, the empty set $\nothing$ is a perfectly good
face of~$\Delta$, representing the empty set of vertices; we set its
dimension equal to~$-1$.  Thus \mbox{$\Delta = \{\nothing\}$} is a
sphere of dimension $-1$, with reduced homology $\ZZ$ in dimension~$-1$.

The notion of vertex-decomposability was introduced by Provan and
Billera \cite{BilleraProvan}, who proved that it implies shellability.
For the convenience of the reader, the proof of the next result is
copied more or less verbatim from \cite[Section~1.8]{grobGeom}.

\begin{thm} \label{thm:shell}
Subword complexes $\Delta(Q,\pi)$ are vertex-decomposable, hence
shellable.
\end{thm}
\begin{proof}
Supposing that $Q = (\sigma, \sigma_2, \sigma_3, \ldots, \sigma_m)$, it
suffices to show that both the link and the deletion of $\sigma$ from
$\Delta(Q,\pi)$ are subword complexes.  By definition, both consist of
subwords of $Q' = (\sigma_2, ..., \sigma_m)$.  The link is naturally
identified with the subword complex $\Delta(Q',\pi)$.  For the deletion,
there are two cases.  If $\sigma \pi$ is longer than~$\pi$, then the
deletion of $\sigma$ equals its link because no reduced expression for
$\pi$ begins with~$\sigma$.  On the other hand, when $\sigma \pi$ is
shorter than~$\pi$, the deletion is $\Delta(Q',\sigma\pi)$.%
\end{proof}

\begin{remark} \label{rk:greedoid}
Among the known vertex decomposable simplicial complexes are the dual
greedoid complexes \cite{BKL}, which include the matroid complexes.
Although subword complexes strongly resemble dual greedoid complexes, the
exchange axioms defining greedoids seem to be slightly stronger than the
exchange axioms for facets of subword complexes imposed by Coxeter
relations.  In particular, the na\"\i ve ways to correspond subword
complexes to dual greedoid complexes do not work, and we conjecture that
they are not in general isomorphic to dual greedoid complexes.

To be precise, a collection $M$ of subsets of a finite vertex set~$Q$
constitutes the \bem{feasible subsets} of a \bem{greedoid} when
$\nothing \in M$, and
\begin{itemize}
\item[]
if $X$ and~$Y$ are in $M$ with $|X| > |Y|$, then there is some element
$x \in X \minus Y$ such that $Y \cup x$ lies in~$M$.
\end{itemize}
The facets of the \bem{dual greedoid complex} are then the complements
in~$Q$ of the maximal elements (\bem{bases}) in~$M$.

There is a natural attempt at defining a greedoid whose dual complex is
$\Delta(Q,\pi)$: namely, let $M(Q,\pi)$ be the collection of subwords
of~$Q$ that are themselves {\em reduced} subwords of some $P \subseteq Q$
representing~$\pi$.  Thus an element $Y \in M(Q,\pi)$ is a sublist of~$Q$
such that (i)~the ordered product of elements in $Y$ has length~$|Y|$,
and (ii)~there is some sublist $Z \subseteq Q$ such that $Y \cup Z$ is a
reduced expression for~$\pi$.  However, this $M(Q,\pi)$ need not be a
greedoid.

An easy non-greedoid example occurs when $\pi = 12543 = s_3s_4s_3 =
s_4s_3s_4$ and $Q$ is the
reduced expression $s_4s_3s_2s_1s_4s_3s_2s_4s_3s_4$ for the long word
in~$S_5$:
$$
\begin{array}{rcc@{}c@{}c@{}c@{}c@{}c@{}c@{}c@{}c@{}c@{}}
        Q &=& s_4 & s_3 & s_2 & s_1 & s_4 & s_3 & s_2 & s_4 & s_3 & s_4 \\
        X &=&     & s_3 &     &     & s_4 & s_3                         \\
        Y &=& s_4 &     &     &     &     &     &     &     & s_3       \\
        Z &=&     &     &     &     &     &     &     &     &     & s_4
\end{array}
$$
Moving any of the elements from~$X$ down to~$Y$ creates a non-reduced
expression.

The reader is invited to find a general construction of greedoids making
subword complexes into dual greedoid complexes; we conjecture that none
exists.  Note that any successful attempt will exclude the subword $Y$
above from the feasible set.
\end{remark}

\section{Balls or spheres}

\label{sec:balls}

Knowing now that subword complexes in Coxeter groups are shellable, we
are able to prove a much more precise statement.  Our proof technique
requires a certain deformation of the group algebra of a Coxeter group.
As we shall see in Lemma~\ref{lemma:dem}.\ref{geq}, the Demazure product
in the following definition ``detects'' Bruhat order on arbitrary words
by a subword condition
just like the ordinary product detects Bruhat order on reduced words by a
subword condition.

\begin{defn} \label{defn:prod}
Let $R$ be a commutative ring, and $\cD$ a free $R$-module with basis
$\{e_{\pi} \mid \pi \in \Pi\}$.  Defining a multiplication on $\cD$ by
\begin{eqnarray} \label{eq:prod}
e_\pi e_\sigma &=& \left\{\begin{array}{@{\,}ll}
        e_{\pi\sigma} & \hbox{if } \length(\pi\sigma) > \length(\pi)\\
        e_{\pi}       & \hbox{if } \length(\pi\sigma) < \length(\pi)
                          \end{array}\right.
\end{eqnarray}
for $\sigma \in \Sigma$ yields the \bem{Demazure algebra} of
$(\Pi,\Sigma)$ over $R$.  Define the \bem{Demazure product} $\delta(Q)$
of the word $Q = (\sigma_1, \ldots, \sigma_m)$ by $e_{\sigma_1} \cdots
e_{\sigma_m} = e_{\delta(Q)}$.
\end{defn}

\begin{example} \label{ex:dem}
Let $\Pi = S_n$
act on the polynomial ring $R[x_1,\ldots,x_n]$ by permuting the
variables.  Define the \bem{Demazure operator $\dem i$} for $i =
1,\ldots,n-1$ on a polynomial $f = f(x_1,\ldots,x_n)$ with coefficients
in~$R$ by
\begin{eqnarray*}
  \dem i(f) = \frac{x_{i+1}f - x_i(s_i\cdot f)}{x_{i+1}-x_i}.
\end{eqnarray*}
Checking monomial by monomial in~$f$ reveals that the denominator divides
the numerator, so this rational function is really a polynomial
in~$R[x_1,\ldots,x_n]$.  The algebra $\cD$ is isomorphic to the algebra
generated over $R$ by the Demazure operators $\dem i$; hence the name
`Demazure algebra'.  In this case, the fact that $\cD$ is an associative
algebra with given free $R$-basis follows from the easily verified fact
that the Demazure operators satisfy the Coxeter relations.
\end{example}

\begin{remark}
The operators in Example~\ref{ex:dem} and the related `divided
difference' operators were introduced by Demazure~\cite{Dem} and
Bernstein--Gel$'$fand--Gel$'$fand \cite{BGG} for arbitrary Weyl groups.
Their context was the calculation of the cohomology and \K-theory
classes of Schubert varieties in $G/P$ via desingularization.  The
operators~$\dem i$, which are frequently denoted in the literature
by~$\pi_i$, were called \bem{isobaric divided differences} by Lascoux and
Sch\"utzenberger \cite{LSgrothVarDrap}.  See Section~\ref{sec:groth} for
the relation to Grothendieck polynomials, and \cite{NoteSchubPoly} for
background on the algebra of divided differences.
\end{remark}

In general, the fact that the equations in~(\ref{eq:prod}) define an
associative algebra is the special case of
\cite[Theorem~7.1]{HumphCoxGrps} where all of the `$a$'~variables
equal~$1$ and all of the `$b$'~variables are zero.  Observe that the
ordered product of a word equals the Demazure product if the word is
reduced.  Here are some basic properties of Demazure products, using
`$\geq$' and `$>$' signs to denote the Bruhat partial order on~$\Pi$, in
which $\tau \geq \pi$ if some (and hence every) reduced word
representing~$\tau$ contains a subword representing~$\pi$
\cite[Section~5.9]{HumphCoxGrps}.  For notation in the proof and
henceforth, we write $Q \minus \sigma_i$ for the word of size $m-1$
obtained from $Q = (\sigma_1,\ldots,\sigma_m)$ by omitting~$\sigma_i$.

\begin{lemma} \label{lemma:dem}
Let $P$ be a word in $\Pi$ and let $\pi \in \Pi$.
\begin{numbered}
\item \label{geq}
The Demazure product $\delta(P)$ is $\geq \pi$ if and only if $P$
contains~$\pi$.

\item \label{=}
If $\delta(P) = \pi$, then every subword of $P$ containing $\pi$ has
Demazure product~$\pi$.

\item \label{>}
If $\delta(P) > \pi$, then $P$ contains a word\/ $T$ representing an
element $\tau > \pi$ satisfying $|T| = \length(\tau) = \length(\pi)+1$.
\end{numbered}
\end{lemma}
\begin{proof}
If $P' \subseteq P$ and $P'$ contains $\pi$, then $P'$ contains
$\delta(P')$ and $\pi = \delta(P) \geq \delta(P') \geq \pi$, proving
part~\ref{=} from part~\ref{geq}.  Choosing any $\tau \in \Pi$ such that
$\length(\tau) = \length(\pi) + 1$ and $\pi < \tau \leq \delta(P)$ proves
part~\ref{>} from part~\ref{geq}.

Now we prove part~\ref{geq}.  Suppose $\pi' = \delta(P) \geq \pi$, and
let $P' \subseteq P$ be the subword obtained by reading $P$ in order,
omitting any reflections along the way that do not increase length.  Then
$P'$ represents $\pi'$ by definition, and contains $\pi$ because any
reduced expression for $\pi'$ contains a reduced expression for $\pi$.

If $T \subseteq P$ represents~$\pi$, then use induction on~$|P|$ as
follows.  Let $\sigma \in \Sigma$ be the last reflection in the list~$P$,
so $\delta(P)\sigma < \delta(P)$ by definition of Demazure product, and
$\delta(P \minus \sigma)$ equals either $\delta(P)$ or $\delta(P)\sigma$.
If $\pi\sigma > \pi$ then $T \subseteq P \minus \sigma$, so $\pi \leq
\delta(P \minus \sigma) \leq \delta(P)$ by induction.  If $\pi\sigma <
\pi$ and $T' \subset T$ represents~$\pi\sigma$, then $T' \subseteq P
\minus \sigma$ and hence $\pi\sigma \leq \delta(P \minus \sigma)$ by
induction.  Since $\pi\sigma < \pi$, we have $\pi\sigma \leq \delta(P
\minus \sigma) \implies \pi \leq \delta(P)$.%
\end{proof}

%
%

\begin{lemma} \label{lemma:T}
Let $T$ be a word in $\Pi$ and $\pi \in \Pi$ such that $|T| =
\length(\pi) + 1$.
\begin{numbered}
\item \label{atmost2}
There are at most two elements $\sigma \in T$ such that $T \minus \sigma$
represents~$\pi$.

\item \label{atleast2}
If\/ $\delta(T) = \pi$, then there are two distinct $\sigma \in T$ such
that $T \minus \sigma$ represents~$\pi$.

\item \label{distinct}
If\/ $T$ represents $\tau > \pi$, then $T \minus \sigma$ represents~$\pi$
for exactly one $\sigma \in T$.
\end{numbered}
\end{lemma}
\begin{proof}
Part~\ref{atmost2} is obvious if $|T| \leq 2$, so suppose there are
elements $\sigma_1,\sigma_2,\sigma_3 \in T$ (in order of appearance) such
that $T \minus \sigma_i$ represents~$\pi$ for each $i = 1,2,3$.  Writing
$T = T_1 \sigma_1 T_2 \sigma_2 T_3 \sigma_3 T_4$, we find that
\begin{eqnarray*}
  T_1 T_2 \sigma_2 T_3 \sigma_3 T_4 &=& T_1 \sigma_1 T_2 T_3 \sigma_3 T_4.
\end{eqnarray*}
Canceling $T_1$ from the left and $T_3 \sigma_3 T_4$ from the right
yields $T_2 \sigma_2 = \sigma_1 T_2$.  It follows that $\pi = T_1
\sigma_1 T_2 \sigma_2 T_3 T_4 = T_1 \sigma_1 \sigma_1 T_2 T_3 T_4 = T_1
T_2 T_3 T_4$, contradicting the hypothesis that $\length(\pi) = |T| - 1$.

In part~\ref{atleast2}, $\delta(T) = \pi$ means there is some $\sigma \in
T$ such that
\begin{romanlist}
\item
$T = T_1 \sigma T_2$;

\item
$T_1 T_2$ represents~$\pi$; and

\item
$\tau_1 > \tau_1 \sigma$, where $T_1$ represents~$\tau_1$.
\end{romanlist}
Omitting some $\sigma'$ from $T_1$ leaves a reduced expression for
$\tau_1 \sigma$ by~(iii).  It follows that $T \minus \sigma'$ and $T
\minus \sigma$ both represent~$\pi$.

Part~\ref{distinct} is the exchange condition.%
\end{proof}

\begin{lemma} \label{lemma:manifold}
Suppose every codimension~1 face of a shellable simplicial
complex~$\Delta$ is contained in at most two facets.  Then $\Delta$ is a
topological manifold-with-boundary that is homeomorphic to either a ball
or a sphere.  The facets of the topological boundary of~$\Delta$ are the
codimension~1 faces of~$\Delta$ contained in exactly one facet
of~$\Delta$.
\end{lemma}
\begin{proof}
\begin{excise}{%
Let $B_d$ and $B_d'$ be homeomorphic to $d$-dimensional balls.  If $h :
B_{d-1} \to B_{d-1}'$ is a homeomorphism between $(d-1)$-balls in their
boundaries, then $B_d \cup_h B_d'$ is again homeomorphic to a $d$-ball.
On the other hand, if $h : \partial B_d \to \partial B_d'$ is a
homeomorphism between their entire boundaries, then $B_d \cup_h B_d'$ is
homeomorphic to a $d$-sphere.  Now induct on the number of facets of
$\Delta$:
if the $d$-sphere case is reached, the assumption on codimension~1 faces
of $\Delta$ implies that the shelling must be complete.%
}\end{excise}%
\cite[Proposition~4.7.22]{BLSWZ}.
\end{proof}



\begin{thm} \label{thm:ball}
The subword complex $\Delta(Q,\pi)$ is a either a ball or a sphere.  A
face $Q \minus P$ is in the boundary of $\Delta(Q,\pi)$ if and only if
$P$ has Demazure product $\delta(P) \neq \pi$.%
\end{thm}
\begin{proof}
That every codimension~1 face of $\Delta(Q,\pi)$ is contained in at most
two facets is the content of part~\ref{atmost2} in
Lemma~\ref{lemma:T}, while shellability is Theorem~\ref{thm:shell}.
This verifies the hypotheses of Lemma~\ref{lemma:manifold} for the first
sentence of the Theorem.

If $Q \minus P$ is a face and $P$ has Demazure product $\neq \pi$, then
$\delta(P) > \pi$ by part~\ref{geq} of Lemma~\ref{lemma:dem}.  Choosing
$T$ as in part~\ref{>} of Lemma~\ref{lemma:dem}, we find by
part~\ref{distinct} of Lemma~\ref{lemma:T} that $Q \minus T$ is a
codimension~1 face contained in exactly one facet of $\Delta(Q,\pi)$.
Thus, using Lemma~\ref{lemma:manifold}, we conclude that $Q \minus P
\subseteq Q \minus T$ is in the boundary of $\Delta(Q,\pi)$.

If $\delta(P) = \pi$, on the other hand, part~\ref{atleast2} of
Lemmas~\ref{lemma:dem} and~\ref{lemma:T} say that every codimension~1
face $Q \minus T \in \Delta(Q,\pi)$ containing $Q \minus P$ is contained
in two facets of $\Delta(Q,\pi)$.  Lemma~\ref{lemma:manifold} says each
such $Q \minus T$ is in the interior of $\Delta(Q,\pi)$, whence $Q \minus
P$~must itself be an interior face.%
\end{proof}

\begin{cor}
The complex $\Delta(Q,\pi)$ is a sphere if $\delta(Q) = \pi$ and a ball
otherwise.
\end{cor}

\section{Hilbert series}

\label{sec:hilb}

Let us review some standard notions from Stanley--Reisner theory.  Fix a
field $\kk$ and a set $\zz = z_1,\ldots,z_m$ of variables.  Suppose
$\Delta$ is a simplicial complex with $m$ vertices, which we think of as
corresponding to the simple reflections $\sigma_1,\ldots,\sigma_m$ in the
word~$Q$.  Recall that the \bem{Stanley--Reisner ideal} of~$\Delta$ is
the ideal $I_\Delta = \<\prod_{i \in D} z_i \mid D \not\in \Delta\>$
generated by monomials corresponding to the (minimal) nonfaces
of~$\Delta$.  Equivalently,
\begin{eqnarray*}
  I_\Delta &=& \bigcap_{D \in \Delta} \<z_i \mid i \not\in D\>
\end{eqnarray*}
is an intersection of prime ideals for faces of~$\Delta$ by an easy
exercise.
By definition, the \bem{Hilbert series} $H(\kk[\Delta];\zz)$ of the
\bem{Stanley--Reisner ring} $\kk[\Delta] = \kk[\zz]/I_\Delta$ equals the
sum of all monomials in $\kk[\zz]$ that lie outside~$I_\Delta$.  Thus
$H(\kk[\Delta];\zz)$ is the sum of all monomials outside every one of the
ideals $\<z_i \mid i \not\in D\>$ for faces $D \in \Delta$.  This sum is
over the monomials $\zz^\bb$ for $\bb \in \NN^m$ having support exactly
$D$ for some face $D \in \Delta$:
\begin{eqnarray} \label{eq:J}
  H(\kk[\Delta];\zz) &=& \sum_{D \in \Delta} \prod_{i \in D}
  \frac{z_i}{1-z_i}
  \ \:=\ \:
  \sum_{D\in\Delta}\frac{\prod_{i\in D}(z_i)\prod_{i\not\in D}(1-z_i)}
  {\prod_{i=1}^m(1-z_i)}. 
\end{eqnarray}

In the special case where $\Delta = \Delta(Q,\pi)$ is a subword complex,
the Stanley--Reisner ideal is
the intersection $I_\Delta = \bigcap \<z_i \mid \sigma_i \in P\>$ over
subwords $P\subseteq Q$ such that $P$ represents~$\pi$.  Now we are ready
to state the main result of this section.

\begin{thm} \label{thm:hilb}
If\/ $\Delta$ is the subword complex $\Delta(Q,\pi)$ and $\ell =
\length(\pi)$, then the Hilbert series of the Stanley--Reisner ring
$\kk[\Delta]$ is
\begin{eqnarray*}
  H(\kk[\Delta]; \zz) &=& \frac{\sum_{\delta(P)=\pi}
  (-1)^{|P| - \ell}(\1-\zz)^P}{\prod_{i=1}^m(1-z_i)},
\end{eqnarray*}
where $(\1-\zz)^P = \prod_{\sigma_i \in P}(1-z_i)$, and the sum is over
subwords $P \subseteq Q$.
\end{thm}

The proof of Theorem~\ref{thm:hilb} is after
Proposition~\ref{prop:inversion}.  First, we set about stating and
proving the two results used in the proof of the theorem.

In general, if $\Gamma$ is an arbitrary monomial ideal $J \subseteq
\kk[\zz]$, or a quotient $\kk[\zz]/J$, then the Hilbert series of
$\Gamma$ (which is the sum of all monomials inside or outside~$J$,
respectively) has the form
\begin{eqnarray*}
  H(\Gamma;\zz) &=& \frac{\KK(\Gamma;\zz)}{\prod_{i=1}^m(1-z_i)},
\end{eqnarray*}
and we call $\KK(\Gamma;\zz)$ the \bem{\K-polynomial} or \bem{Hilbert
numerator} of~$\Gamma$.  It has the following direct interpretation in
terms of $\ZZ^m$-graded homological algebra.  Since $\Gamma$ is
$\ZZ^m$-graded, it has a minimal $\ZZ^m$-graded free resolution
\begin{equation} \label{eq:freeres}
  0 \from \Gamma \from E_0 \from E_1 \from\cdots\from E_m \from 0,\qquad
  E_j = \bigoplus_{P \subseteq Q} \kk[\zz](-\deg \zz^P)^{\beta_{j,P}},
\end{equation}
where $\beta_{j,P}$ is the \bem{$j^\th$ Betti number} of $\Gamma$ in
$\ZZ^m$-graded degree $\deg\zz^P$.  Then the \K-polynomial of $\Gamma$
is $\KK(\Gamma; \zz) = \sum_j(-1)^j \beta_{j,P} \cdot \zz^P$.

Hochster's formula, which shall state in~\eqref{eq:hoc}, says how to
calculate explicitly the Betti numbers of the \bem{Alexander dual ideal}
of~$I_\Delta$, which is defined by
\begin{eqnarray*}
  I_\Delta^\star &=& \Big\<\prod_{i \not\in D} z_i \mid D \in
  \Delta\Big\>.
\end{eqnarray*}
Note that the generators of $I_\Delta^\star$ are obtained by multiplying
the variables in each prime component of~$I_\Delta$.  Thus, for instance,
when $\Delta = \Delta(Q,\pi)$ is a subword complex,~we~get
\begin{eqnarray*}
  I_\Delta^\star &=& \<\zz^P \mid P\subseteq Q \hbox{ and } P \hbox{
  represents } \pi\>,
\end{eqnarray*}
where $\zz^P = \prod_{\sigma_i \in P} z_i$ for any subword $P \subseteq
Q$.  Hochster's formula \cite[p.~45]{MP} now says that, in terms of
reduced homology of $\Delta = \Delta(Q,\pi)$ over the field~$\kk$, the
$\ZZ^m$-graded Betti numbers of $I_\Delta^\star$ over $\kk[\zz]$ are
\begin{eqnarray} \label{eq:hoc}
  \beta_{j,P} &=& \dim_\kk\HH_{j-1}(\link(Q \minus P,\,\Delta);\kk).
\end{eqnarray}

\begin{lemma} \label{lemma:links}
If\/ $\Delta$ is the subword complex $\Delta(Q,\pi)$ and $\ell =
\length(\pi)$, then
\begin{eqnarray*}
  \KK(I_\Delta^\star; \zz) &=& \sum_\twoline{P \subseteq Q}
  {\delta(P)=\pi} (-1)^{|P|-\ell} \zz^P
\end{eqnarray*}
is the Hilbert numerator of the Alexander dual ideal.
\end{lemma}
\begin{proof}
Let $Q \minus P \in \Delta$, so $P \subseteq Q$ contains~$\pi$.  By
Theorem~\ref{thm:ball}, either $\delta(P) \neq \pi$, in which case
$\link(Q \minus P,\Delta)$ is contractible, or $\delta(P) = \pi$, in
which case $\link(Q \minus P,\Delta)$ is a sphere of dimension
\begin{eqnarray*}
  (\dim \Delta) - |Q \minus P| &=& (|Q| - \ell - 1) - |Q\minus P|\
  \:=\ \: |P| - \ell - 1.
\end{eqnarray*}
(Recall that a sphere of dimension~$-1$ is taken to mean the empty
complex~$\{\nothing\}$ having nonzero reduced homology in
dimension~$-1$.)  Therefore $\HH_{j-1}\link(Q \minus P,\Delta)$ is zero
unless $\delta(P) = \pi$ and $j = |P| - \ell$, in which case the reduced
homology has dimension~$1$.  Now apply (\ref{eq:hoc}) to the formula
$\KK(\Gamma; \zz) = \sum_j(-1)^j \beta_{j,P} \cdot \zz^P$.%
\end{proof}

%

Lemma~\ref{lemma:links} helps us calculate the Hilbert series
of~$\kk[\Delta]$ because the \K-polynomials of the Stanley--Reisner
ring $\kk[\Delta]$ and the Alexander dual ideal $I_\Delta^\star$ are
intimately related, as the next result demonstrates.  Although it holds
more generally for the ``squarefree modules'' of Yanagawa \cite{Yan}, as
shown in \cite[Theorem~4.36]{Mil3}, we include an elementary proof of
Proposition~\ref{prop:inversion} because of its simplicity.
A~$\ZZ$-graded version
was proved by Terai for squarefree ideals using some calculations
involving $f$-vectors of simplicial complexes \cite[Lemma~2.3]{Ter}.
For notation, $\KK(\Gamma;\1-\zz)$ is the polynomial obtained from
$\KK(\Gamma;\zz)$ by substituting $1-z_i$ for each variable~$z_i$.

\begin{prop}[Alexander inversion formula] \label{prop:inversion}
For any simplicial complex $\Delta$ we have $\KK(\kk[\Delta];\zz) =
\KK(I_\Delta^\star;\1-\zz)$.
\end{prop}
\begin{proof}
See~(\ref{eq:J}) for the Hilbert series of $\kk[\Delta]$.  On the other
hand, the Hilbert series of $I_\Delta^\star$ is the sum of all monomials
$\zz^\bb$ divisible by $\prod_{i \not\in D} z_i$ for some $D \in \Delta$:
\begin{equation} \label{eq:Jstar}
  H(I_\Delta^\star;\zz)\ \:=\ \:\sum_{D \in \Delta} \prod_{i \not\in D}
  \frac{z_i}{1-z_i}
  \ \:=\ \:
  \sum_{D \in \Delta}\frac{\prod_{i\not\in D}(z_i)\prod_{i\in D}(1-z_i)}
  {\prod_{i=1}^m(1-z_i)}.
\end{equation}
Now compare the last expressions of (\ref{eq:J}) and~(\ref{eq:Jstar}).
\end{proof}

\begin{proofof}{Theorem~\ref{thm:hilb}}
Lemma~\ref{lemma:links} gives the numerator of the Hilbert series of the
Alexander dual ideal~$I_\Delta^\star$ (to be defined, below), and
Proposition~\ref{prop:inversion} says how to recover the numerator of the
Hilbert series of the Stanley--Reisner ring from that.
\end{proofof}



Theorem~\ref{thm:hilb} can be restated in a somewhat different form,
grouping subwords with Demazure product~$\pi$ according to their
lexicographically first reduced subwords for~$\pi$.  Given a reduced
subword $D \subseteq Q$, say that $\sigma_i \in Q \minus D$ is
\bem{absorbable} if the word $T = D \cup \sigma_i$ in~$Q$ has the
properties: (i)~\mbox{$\delta(T) = \delta(D)$},
and (ii)~the unique reflection $\sigma_j \in D$ (afforded by
Lemma~\ref{lemma:T}.\ref{atleast2}) satisfying $\delta(T \minus \sigma_j)
= \delta(D)$ has \mbox{index $j < i$}.

\begin{thm} \label{thm:hilb'}
If $\Delta$ is the subword complex $\Delta(Q,\pi)$ and\/ $\abs(D)
\subseteq Q$ is the set of absorbable reflections for each reduced
subword \mbox{$D \subseteq Q$}, then\/ $\kk[\Delta]$ has \K-polynomial
\begin{eqnarray*}
  \KK(\kk[\Delta]; \zz) &=& \sum_{\mathrm{facets\ }Q \minus D} (\1-\zz)^D
  \zz^{\abs(D)},
\end{eqnarray*}
where $(\1-\zz)^D = \prod_{\sigma_i \in D}(1-z_i)$, and\/ $\zz^{\abs(D)}
= \prod_{\sigma_i \in \abs(D)} z_i$.
\end{thm}
\begin{proof}
Given a subword $P \subseteq Q$, say that $P$ \bem{simplifies} to~$D
\subseteq P$, and write \mbox{$P\goesto D$}, if $D$ is the
lexico\-graphically first subword of $P$ with Demazure product
$\delta(P)$.  If~$P$ has Demazure product~$\pi$ and $P \goesto D$, the
subword $Q \minus D$ is automatically \mbox{a facet of~$\Delta$}.

If we denote by $P_{\leq i}$ the initial string of reflections in
$P$ with index at most~$i$, the simplification $D$ is obtained from $P$
by omitting any reflection $\sigma_i \in P$ such that
$\delta(P_{\leq i-1}) = \delta(P_{\leq i})$.  
Theorem~\ref{thm:hilb} says that
\begin{eqnarray*}
  \KK(\kk[\Delta]; \zz) &=& \sum_{\mathrm{facets\ }Q \minus D} (\1-\zz)^D
  \sum_{P \goesto D} (\zz-\1)^{P \minus D}.
\end{eqnarray*}
Now note that subwords $P$ simplifying to~$D$ are (by definition of
Demazure product) obtained by adding to~$D$ (at will) some of its
absorbable reflections in~$Q$.  Therefore
$$
  \sum_{P \goesto D} (\zz-\1)^{P \minus D}\:\ =\:\ \prod_{\sigma_i \in
  \abs(D)} \big(1 + (z_i - 1)\big)
  \:\ =\:\  \zz^{\abs(D)},
$$
completing the proof.
\end{proof}

\begin{remark}
The Hilbert numerator as expressed in Theorem~\ref{thm:hilb'} looks more
like one would expect from a shellable simplicial complex, using a
version of \cite[Proposition~2.3]{Sta} suitably enhanced for the fine
grading.  We believe the reason comes from the \bem{facet adjacency
graph}~$\Gamma(Q,\pi)$ of $\Delta(Q,\pi)$, which by definition has the
facets of~$\Delta(Q,\pi)$ for vertices, while its edges are the interior
\bem{ridges} (codimension~$1$ faces) of~$\Delta(Q,\pi)$.  Two facets are
adjacent if they share a ridge.  Note that every interior ridge lies in
exactly two facets by Lemma~\ref{lemma:manifold}.

The facet adjacency graph $\Gamma(Q,\pi)$ can be oriented, by having
each ridge $Q \minus P$ point toward the facet $Q \minus D$ whenever $P$
simplifies to~$D$.  The resulting directed facet adjacency graph is
acyclic~--- so its transitive closure is a poset~--- because the
relation by ridges is a subrelation of lexicographic order.  We believe
that every linear extension of this poset gives a shelling order for
$\Delta(Q,\pi)$.  The shelling formulae we get for the \K-polynomial
will all be the same, namely the one in Theorem~\ref{thm:hilb'}.
\end{remark}

\begin{excise}{
  
  \begin{prop}\label{prop:shelling}
  The directed facet adjacency graph~$\Gamma(Q,\pi)$ is acyclic, so its
  transitive closure is a poset.  Every linear extension of this poset
  gives a shelling of $\Delta(Q,\pi)$,
  by listing facets in increasing order.
  \end{prop}
  
  Of course, we already knew these complexes to be shellable (Theorem
  \ref{thm:shell}).

  \begin{proof}
  $\Gamma(Q,\pi)$ is acyclic because the relation by ridges is a
  subrelation of lexicographic order.  Let~$\prec$ be an arbitrary
  linear extension.
  
  What we need: If $Q \minus D \prec Q \minus E$ are facets, and $P = D
  \cup E$, so that $Q \minus P = (Q \minus D) \cap (Q \minus E)$, then
  we must be able to find some other facet $Q \minus F \prec Q \minus E$
  such that $F \cup E \subseteq P$ and $|F \cup E| = \length(\pi)+1$.
  
  Let $Q\minus P$ be an interior face of $\Delta(Q,\pi)$, so~$P$ has
  Demazure product~$\pi$, and let~$D$ be the reduced subword to which
  $P$ simplifies.  The face $Q\minus P$ might lie inside many facets
  ($P$~might contain many reduced expressions), but we claim now that
  $Q\minus D$ is the $\prec$-least with product~$\pi$.
  
  Let $D' \neq D$ be a reduced subword of~$P$ with product~$\pi$, and
  $\sigma_d \in D \minus D'$ the first letter of~$D$ not in~$D'$.  Set
  $E = D' \cup \sigma_d$.  Since $D \subset E \subseteq P$, the Demazure
  product of~$E$ is again~$\pi$.  Hence by
  Lemma~\ref{lemma:T}.\ref{atleast2} there is a unique other reflection
  $\sigma_{d'}$ in~$D'$ satisfying $\delta(E \minus \sigma_{d'}) = \pi$,
  and $d' > d$ because $D$ is the simplification of~$P$.  It follows
  that $Q\minus E$ is an edge of $\Gamma(Q,\pi)$ connecting the vertex
  $Q\minus D'$ in~$\Gamma(Q,\pi)$ to the vertex~$Q\minus (E \minus d')$,
  and oriented away from $Q\minus D'$.  We conclude that $D \prec D'$ by
  induction on the length of the initial segements shared by~$D$
  and~$D'$.
  
  In particular, this tells us at which point $Q\minus P$ is added to
  the growing complex during the purported shelling that $\prec$ gives
  us: it's when we add the facet $Q\minus D$, $P\goesto D$.  So when we
  attach the facet $Q\minus D'$, its intersection with the previously
  added $Q\minus D$ is contained inside the ridge $Q\minus E$.
  
  More generally, if $Q\minus P$ is a boundary face, let $P' \subseteq
  P$ be $D$ union the {\em backwards-absorbable} letters $\sigma_i$ --
  the same definition as absorbable, except with $j>i$ rather than
  $j<i$.  Then $Q\minus P'$ is an interior face, contained in $Q\minus
  D$, and containing $Q\minus P$.
  
  At this point we need to know that if $D$ is not the
  lex-first $\pi$-word in $P$, then $P'$ is strictly bigger than
  $D$. Why is that?
  \end{proof}
  
  So while every choice of a shelling~$\prec$ gives us in a standard way
  a formula for the \K-polynomial of the complex, we always get the same
  formula
  
}\end{excise}%

\section{Combinatorics of Grothendieck polynomials}

\label{sec:groth}

The Grothendieck polynomial \hbox{$\GG_w(\xx)$} in variables
$x_1,\ldots,x_n$ and its ``double'' analogue $\GG_w(\xx,\yy)$ represent
the classes of Schubert varieties in ordinary and equivariant \K-theory
of the flag manifold \cite{LSgrothVarDrap}.  Their algebraic definition
will be recalled below.  The goal of this section is to derive as special
cases of Theorem~\ref{thm:hilb} and \ref{thm:hilb'} two formulae for
Grothendieck polynomials.  The first formula (Corollary~\ref{cor:links})
coincides with a special case of a formula discovered by Fomin and
Kirillov \cite{FKgrothYangBax}.  It is the \K-theoretic analogue of the
Billey--Jockusch--Stanley formula for the Schubert
polynomial~$\SS_w(\xx)$ \cite{BJS,FSnilCoxeter}, interpreted here for the
first time in terms of simplicial topology.  The second formula
(Corollary~\ref{cor:links'}) relates to other combinatorial models for
Grothendieck polynomials in work of Lenart, Robinson, and Sottile
\cite{LenRobSot}.  We begin with the example of subword complexes that
pervades \cite{grobGeom}.

\begin{example} \label{ex:square}
Set $\Pi = S_{2n}$, and let $Q_{n \times n} =$
$$
  s_ns_{n-1} \ldots s_2s_1\ s_{n+1}s_n \ldots s_3s_2\
  s_{n+2}s_{n+1}\ldots \ \ \ldots s_{n+2}s_{n+1}\ s_{2n-1}s_{2n-2} \ldots
  s_{n+1}s_n.
$$
This is the \bem{square word} from \cite[Example~1.8.3]{grobGeom}, so
named because the $n^2$ simple reflections in this list~$Q$ fill the $n
\times n$ grid naturally by starting at the upper-right, continuing to
the left, and subsequently reading each row from right to left, in turn.
Observe that every occurrence of $s_i$ in $Q_{n \times n}$ sits on the
$i^\th$ antidiagonal of the resulting square array.

Given $w \in S_n$ (not $S_{2n}$),
the subword complex $\Delta = \Delta(Q_{n \times n},w)$ plays a crucial
role in the main theorems of \cite{grobGeom}; see
Proposition~\ref{prop:square}, below.  For the Stanley--Reisner ring
$\kk[\Delta]$, we index the variables $\zz = z_1,\ldots,z_{n^2}$ by
their positions in the $n \times n$ grid, so $\zz =
\{z_{ij}\}_{i,j=1}^n$ with $z_{11}$ at the upper-left and $z_{1n}$ at
the upper-right.%
\end{example}


\begin{defn} \label{defn:groth}
Let $w \in S_n$ be a permutation, and recall the Demazure operators $\dem
i$ from Example~\ref{ex:dem}.  The \bem{Grothendieck polynomial}
$\GG_w(\xx)$ is obtained recursively from the top one $\GG_{w_0}(\xx) :=
\prod_{i=1}^n (1-x_i)^{n-i}$ via the recurrence
\begin{eqnarray*}
  \GG_{ws_i}(\xx) &=& \dem i \GG_w(\xx)
\end{eqnarray*}
whenever $\length(ws_i) < \length(w)$.  The \bem{double Grothendieck
polynomials} are defined by the same recurrence, but start from
$\GG_{w_0}(\xx,\yy) := \prod_{i+j \leq n} (1-x_i y_j)$.
\end{defn}

We use slightly different notation in Definition~\ref{defn:groth} than
in \cite[Definition~1.1.3]{grobGeom}: the polynomial $\GG_w(\xx,\yy)$
here is obtained from the corresponding Laurent polynomial in
\cite{grobGeom} by setting each variable $y_i^{-1}$ to~$y_i$ (the
geometry in \cite{grobGeom} required inverses).  This alteration makes
the notation more closely resemble that in \cite{FKgrothYangBax}, where
their polynomial $\mathfrak L\begin{array}{@{}l@{}}\\[-4ex]\scriptstyle
(-1)\\[-1.5ex] \scriptstyle\, w\\[-1ex]\end{array}(y,x)$ corresponds to
what we call here $\GG_w(\1-\xx,\1-\yy)$.

Grothendieck polynomials connect to subword complexes by part of the
`Gr\"obner geometry theorems' in \cite{grobGeom}.  In our context, it
says the following.

\begin{prop}[{\cite[Theorem~A]{grobGeom}}] \label{prop:square}
Suppose $w \in S_n$, and let $\Delta = \Delta(Q_{n \times n},w)$ be the
subword complex for the square word.  Setting $z_{ij}$ equal to $x_iy_j$
or to~$x_i$ in the Hilbert numerator $\KK(\kk[\Delta];\zz)$ yields
respectively the double Grothendieck polynomial $\GG_w(\xx,\yy)$ or the
Grothendieck polynomial~$\GG_w(\xx)$.
\end{prop}

For notation, regard subwords $P \subseteq Q_{n \times n}$ as subsets of
the $n \times n$ grid.

\begin{cor}[{\cite[Theorem~2.3 and p.~190]{FKgrothYangBax}}]\label{cor:links}
If $w \in S_n$ and $Q_{n \times n}$ is the square word as in
Example~\ref{ex:square}, then the double Grothendieck polynomial
$\GG_w(\xx,\yy)$ satisfies
\begin{eqnarray*}
  \GG_w(\1-\xx,\1-\yy) &=& \sum_\twoline{P \subseteq Q_{n \times
  n}}{\delta(P) = w}\prod_{(i,j)\in P} (-1)^{|P| - \ell} (x_i + y_j - x_i
  y_j),
\end{eqnarray*}
where\/ $\length(w) = \ell$.  The version for single Grothendieck
polynomials reads
\begin{eqnarray*}
  \GG_w(\1-\xx) &=&
  \sum_\twoline{P \subseteq Q_{n \times n}}{\delta(P)=w} (-1)^{|P| -
  \ell} \xx^P, \qquad\hbox{where}\qquad \xx^P = \prod_{(i,j) \in P} x_i.
\end{eqnarray*}
\end{cor}
\begin{proof}
Apply Theorem~\ref{thm:hilb} to the subword complex $\Delta = \Delta(Q_{n
\times n},w)$.  Substituting $x_i y_j$ for $z_{ij}$ as stipulated in
Proposition~\ref{prop:square} yields the double version after calculating
$1-(1-x_i)(1-y_j) = x_i+y_j-x_i y_j$, while the single version follows
trivially.%
\end{proof}


Corollary~\ref{cor:links} can be rewritten in terms of
absorbable reflections as in Theorem~\ref{thm:hilb'}.

\begin{cor} \label{cor:links'}
If $w \in S_n$ and $\Delta(Q_{n \times n},w)$ is the square subword
complex, then
\begin{eqnarray*}
  \GG_w(\xx,\yy) &=& \sum_\twoline{\mathrm{facets\ }}{Q_{n\times n}
  \minus D} \,\prod_{\:(i,j)\in D} (1 - x_i y_j) \prod_{(i,j) \in
  \abs(D)} x_i y_j,
\end{eqnarray*}
The version for single Grothendieck polynomials reads
\begin{eqnarray*}
  \GG_w(\xx) &=& \sum_\twoline{\mathrm{facets\ }}{Q_{n\times n} \minus D}
  (\1-\xx)^D \xx^{\abs(D)},
\end{eqnarray*}
where $(\1-\xx)^D = \prod_{(i,j) \in D} (1-x_i)$ and\/ $\xx^{\abs(D)} =
\prod_{(i,j) \in \abs(D)} x_i$.
\end{cor}
\begin{proof}
Apply Proposition~\ref{prop:square} to the result of
Theorem~\ref{thm:hilb'} for $\Delta = \Delta(Q_{n \times n},w)$.%
\end{proof}

Readers familiar with reduced pipe dreams (also called rc-graphs; see
\cite[Section~1.4]{grobGeom} for an introduction) can see a geometric
interpretation of
absorbable reflections: given a reduced pipe dream~$D$, an elbow tile is
absorbable if the two pipes passing through it intersect in a crossing
tile to its northeast.  Thus the pipe dream~$D$ ``nearly misses'' being a
reduced pipe dream because of that elbow tile.  Note that there must be
exactly one reduced pipe dream with no
absorbable elbow tiles (in \cite{BB} this is the `bottom' rc-graph),
because the constant term of the \K-polynomial of any Stanley--Reisner
ring~--- or indeed any quotient of the polynomial ring~---~equals~$1$.

Here is a weird consequence of the Demazure product characterization of
the Hilbert numerator $\KK(\kk[\Delta];\zz)$ for $\Delta = \Delta(Q_{n
\times n},w)$

\begin{por} \label{por:weird}
For each squarefree monomial $\zz^P$ in the variables $\zz =
(z_{ij})_{i,j=1}^\infty$, there exists a unique permutation $w \in
S_\infty = \bigcup_n S_n$ such that $\zz^P$ appears with nonzero
coefficient in the Hilbert numerator of the Alexander dual ideal
$I_{\Delta(Q_{n \times n},w)}^\star$ for some (and hence all) $n$ such
that $w \in S_n$.  The coefficient of $\zz^P$ is $\pm 1$.
\end{por}
\begin{proof}
The permutation $w$ in question is $\delta(P)$, by
Lemma~\ref{lemma:links}.
\end{proof}

\section{Open problems}

\label{sec:open}

The considerations in this paper motivate some questions concerning the
combinatorics of reduced expressions in Coxeter groups.

\begin{question} \label{q:1}
Given an element $\pi \in \Pi$, what is the smallest size of a word in
$\Sigma$ containing every reduced expression for $\pi$ as a subword?
\end{question}
\noindent
Note that a smallest size word
containing all reduced expressions for~$\pi$ will not in general be
unique.  Indeed, even for the long word $w = 321 \in S_3$, there are two
such: $s_1s_2s_1s_2$ and $s_2s_1s_2s_1$.

Question~\ref{q:1} asks for a measure of how far intervals in the weak
order are from being subword complexes.  Another measure would be
provided by a solution to the following problem, which asks roughly how
many faces must be added to order complexes of intervals in the weak
order to get subword complexes.  To be precise, let the \bem{repetition
number} $\repnum(Q,\pi)$ be the largest number of times that a single
reduced expression for~$\pi$ appears as a subword of~$Q$.

\begin{problem} \label{q:2}
Describe the function sending $\pi \mapsto \repnum(\pi)$, where
$\repnum(\pi) = \min(\repnum(Q) \mid Q$ contains all reduced expressions
for $\pi)$.
\end{problem}

Restricting to symmetric groups, for instance,
\begin{question} \label{q:3}
Is the function in Problem~\ref{q:2} bounded above on $S_\infty =
\bigcup_n S_n$?  If not, how does it grow?
\end{question}

Given that subword complexes appeared naturally in the context of the
geometry of Schubert varieties, it is natural to ask whether there are
good geometric representatives for subword complexes.

\begin{question} \label{q:4}
Can any spherical subword complex be realized as a convex polytope?
\end{question}

One could also take the opposite perspective, by starting with a
simplicial sphere.

\begin{problem} \label{q:5}
Characterize those simplicial spheres realizable as subword complexes.
\end{problem}

Of course, in all of these problems
it may be useful to try restricting to words in~$S_n$.

\vspace{-.1ex}

\def\cprime{$'$}
\providecommand{\bysame}{\leavevmode\hbox to3em{\hrulefill}\thinspace}
\raggedbottom


\end{document}